%dt1.tex: Boolean  Function Analogs of Covering Systems

%begin macros

\baselineskip=14pt
\parskip=10pt

\magnification=\magstephalf

\def\1{{\overline{1}}}
\def\2{{\overline{2}}}
\parindent=0pt
\overfullrule=0in

\def\frac#1#2{{#1 \over #2}}
%\headline={\rm  \ifodd\pageno  \RightHead  \else  \LeftHead  \fi}
%\def\RightHead{\centerline{
%Title
%}}
%\def\LeftHead{ \centerline{Doron Zeilberger}}
%end macros
\centerline
{\bf Boolean  Function Analogs of Covering Systems}
\bigskip
\centerline
{\it Anthony ZALESKI and Doron ZEILBERGER}

{\bf Abstract}: Bob Hough recently disproved a long-standing conjecture of Paul Erd\H{o}s regarding
covering systems. Inspired by his seminal paper, we describe analogs of covering systems to
Boolean functions, and more generally, the problem of covering discrete 
hyper-boxes by non-parallel lower dimensional hyper-subboxes. We point out that
very often primes are red herrings. This is definitely the case for covering system, and
who knows, perhaps also for the Riemann Hypothesis.

\bigskip

{\bf Prime Numbers are Sometimes Red Herrings}

The great French mathematical columnist Jean-Paul Delahaye [D] recently posed
the following brain-teaser, adapting a beautiful puzzle, of unknown origin,
popularized by  Peter Winkler in his wonderful book [W] (pp. 35-43).

Here is a free translation from the French.

{\bf Enigma: Nine Beetles and prime numbers}

One places nine beetles on a circular track, where the nine arc distances,
measured in meters, between two consecutive beetles are the first nine prime numbers,
2,3,5,7,11,13,17,19 and 23. The order is arbitrary, and each number appears exactly
once as a distance.

At starting time, each beetle decides {\it randomly} whether she would go, traveling at a speed of $1$ meter
per minute, clockwise or counter-clockwise. When two beetles bump into each other, they immediately
do a ``U-turn," i.e. reverse direction. We assume that the size of the beetles
is negligible. At the end of $50$ minutes, after many collisions, one notices the
distances between the new positions of the beetles. The nine distances are exactly as before,
the first nine prime numbers! How to explain this miracle?

Before going on to the next section, we invite you to solve this lovely puzzle all by yourself.

{\bf Solution of the Enigma}

Note that the length of the circular track is $2+3+5+7+11+13+17+19+23 \,= \, 100$ meters.

Let each beetle carry a flag, and whenever they bump into each other, let them exchange flags.
Since the flags always move in the same direction, and  also move at a speed of $1$ meter per minute, after
$50$ minutes, each flag is {\bf exactly} at the ``antipode" of its original location; hence, the distances are the same!
Of course, this works if the original distances were {\it any}  sequence of numbers: All that they have
to obey is that their sum equals $100$, or more generally, that 
half the sum of the distances divides the product of the speed ($1$ meter per minute in this puzzle) and the elapsed time ($50$ minutes in this puzzle).

This variation, due to Delahaye, is {\bf much} harder than the original version 
posed in [W], where also the initial distances were arbitrary.
In Delahaye's rendition, the solver is bluffed into trying to use the fact that the distances
are primes. Something analogous happened to the great Paul Erd\H{o}s, the patron saint
of combinatorics and number theory, who introduced {\bf covering systems}.

{\bf Covering Systems}

In 1950, Paul Erd\H{o}s [E1], introduced the notion of {\bf covering systems}. A {\it covering system} is
a finite set of arithmetical progressions
$$
\{ \, a_i \, (\,\, mod \, m_i) \quad \vert \quad 1 \leq i \leq N \, \} \quad ,
$$
whose union is the set of all non-negative integers. For example
$$
\{ 0 \, (\, \, mod \, 1) \} \quad ,
$$
is such a (not very interesting) covering system, while
$$
\{ \, 0\, ( \, \,mod \, 2) \quad , \quad  1\, ( \, \, mod \, 2) \}\quad ,
$$
and
$$
\{\, 0 \, ( \, mod \, 5) \quad , \quad  1\, ( \, mod \, 5) \quad , \quad 2 \, ( \, mod \, 5) \quad , \quad 3 \, ( \, mod \, 5) \quad , 
\quad 4 \, ( \, mod \, 5) \, \} \quad ,
$$
are other, almost as boring examples. A slightly more interesting example is
$$
\{ \, 0 \, ( \, mod \, 2) \quad, \quad 1 \, ( \, mod \, 4) \quad, \quad  3 \, ( \, mod \, 4) \, \} \quad .
$$

A covering system is {\bf exact} if all the congruences are disjoint (like in the above boring examples). It is {\bf distinct} if
all the moduli are different.   [From now on, let $a\,(\,b\,)$ mean $a \, ( \, mod \, b)$.]

Erd\H{o}s gave the smallest possible example of a distinct covering system:
$$
\{ \, 0\,(\,2\,) \quad, \quad 0\,(\,3\,) \quad, \quad 1\,(\,4\,) \quad, \quad 5\,(\,6\,)  \quad, \quad 7\,(\,12\,) \, \} \quad .
$$

Of course, the above covering system is not exact since, for example, $0\,(\,2\,)$ and $0\,(\,3\,)$ both contain any multiple of $6$.
A theorem proved by Mirsky and (Donald) Newman, and independently by Davenport and Rado (described in [E2])
implies that a covering system cannot be both exact {\bf and} distinct.
Even a stronger statement holds. Assuming that
our system $\{a_i( m_i)\}_{i=1}^N$ is written in non-decreasing order of the moduli
$m_1 \leq m_2 \leq \dots \leq m_N$, the Mirsky-Newman-Davenport-Rado theorem asserts that $m_{N-1}=m_N$, in other words,
the two top moduli are equal (and hence an exact covering system can never be distinct).
See [Zei1] for an exposition of their snappy proof. While their proof was nice, it was
not as nice as the combinatorial-geometrical proof
that was found by Berger, Felzenbaum, Fraenkel ([BFF1][BFF2]), and exposited in [Zei1].
In fact, they proved the more general Znam theorem that asserts that the highest moduli
shows up at least $p$ times, where $p$ is the smallest prime dividing $lcm\,(\, m_1, \dots, m_N\,)$.
Jamie Simpson ([S]) independently  found a similar proof, but unfortunately chose not to express it in
the evocative geometrical language.

\vfill\eject

{\bf The Berger-Felzenbaum revolution: From Number Theory to Discrete Geometry via the Chinese Remainder Theorem}

While it is a sad truth that the set of positive integers is an {\it infinite} set, a covering system is a 
{\bf finite} object. In order to verify that a  covering system, $\{a_i( m_i)\}_{i=1}^N$ is indeed one,
it suffices to check that it covers all the integers $n$ between $0$ and $M-1$, where
$$
M \,  = \, lcm \, (\, m_1, m_2, \dots , m_N \, ) \quad .
$$

By the fundamental theorem of arithmetic
$$
M=p_1^{r_1} p_2^{r_2} \cdots p_k^{r_k} \quad ,
$$
where $p_1, \dots, p_k$ are primes and $r_1, \dots, r_k$ are positive integers.

For the sake of simplicity, let's assume that $M$ is square-free, i.e. all the exponents $r_1, \dots, r_k$ equal $1$.
The same reasoning, only slightly more complicated, applies in the general case.

Now we have
$$
M=p_1 p_2 \cdots p_k \quad .
$$

The ancient, but still useful, {\bf Chinese Remainder Theorem} tells you that there is a bijection between
the set of integers between $0$ and $M-1$, let's call it $[0,M-1]$, and the Cartesian product of $[0,p_i-1]$, $i=1 \dots k$.
$$
f:= [0, M-1] \rightarrow \prod_{i=1}^{k} [0, p_i -1] \quad,
$$
defined by
$$
f(x):= \, [x \,\, mod \,\, p_1 \, ,  \, x \,\, mod \,\, p_2 \, , \, \quad \dots \quad ,x \,\, mod \,\, p_k \, ] \quad .
$$
So each integer in $[0, M-1]$ is represented by a {\bf point} in the $p_1 \times p_2 \times \dots \times p_k$ 
$k$-dimensional discrete box $ \prod_{i=1}^{k} [0, p_i -1]$.

If $a\,(\,m\,)$ is a member of our covering system, since $m$ is a divisor of $M$, it can be written
as a product of some of the primes in $\{p_1, \dots, p_k\}$, say
$$
m \,= \, p_{i_1} \,p_{i_2} \, \cdots \, p_{i_s} \quad .
$$
Let
$$
m_{i_1}\, = \, a \, mod \,\, p_{i_1} \quad , \quad
m_{i_2} \,=\,  a \, mod \, \, p_{i_2} \quad , \quad
\dots \quad, \quad
m_{i_s} \, = \, a \, mod \,\,  p_{i_s} \quad .
$$
It follows that the members of the congruence $a(m)$ correspond to the points in the $k-s$-dimensional {\bf subbox}
$$
\{(x_1, \dots, x_k) \in [0,p_1-1] \times \dots \times [0,p_k-1] \quad \vert \quad x_{i_1}=m_{i_1}, \quad \dots \quad , x_{i_s}=m_{i_s} \} \quad.
$$
For example if $M=30=2 \cdot 3 \cdot 5$, the congruence class $7(10)$, corresponds to the one-dimensional subbox
(since $ 7 \, mod \,2 =1$ and $7 \, mod \, 5 =2$)
$$
\{ \, (x_1,x_2,x_3) \, \vert \, x_1=1 \quad, \quad 0 \leq x_2 \leq 2 \quad , \quad x_3=2 \} \quad .
$$

In other words a covering system (with square-free $M$) is {\bf nothing but} a way of expressing a certain
$k$-dimensional discrete box as a union of sub-boxes. This was the beautiful insight of Marc Berger,
Alex Felzenbaum, and Aviezri Fraenkel, nicely exposited in [Zeil1].

{\bf Erd\H{o}s's Famous Problem and Bob Hough's Refutation}

 Erd\H{o}s ([E2]) famously asked whether there exists a distinct covering system
$$
a_i (\,\, mod \, m_i) \quad , \quad 1 \leq i \leq N \quad , \quad m_1 <m_2 < \dots <m_N \quad,
$$
with the smallest modulo, $m_1$, arbitrarily large.

As computers got bigger and faster, people (and their computers) came up with examples
that progressively made $m_1$ larger and larger, and many humans thoughts that indeed $m_1$ can be made
as large as one wishes. This was brilliantly refuted by Bob Hough ([H]) who proved that
$m_1 \leq 10^{16}$. This is definitely not sharp, and the true largest $m_1$ is probably less than
$1000$.

Let's now move on from  number theory to something apparently very different: logic!

{\bf Boolean Functions}

Let's recall some basic definitions. 
A {\bf Boolean function} (named after George Boole ([Bo])) of $n$ variables is a function from $\{False,True\}^n$ to $\{False,True\}$.
Altogether there are $2^{2^n}$ Boolean functions of $n$  variables. Any Boolean function $f(x_1, \dots, x_n)$, is determined by  its {\it truth table},
or equivalently, by the set $f^{-1}(True)$, one of the $2^{2^n}$ subsets of   $\{False,True\}^n$.

The {\it simplest} Boolean functions are the
{\bf constant} Boolean functions {\bf True} (the {\it tautology}) corresponding to the whole of $\{False,True\}^n$, and
{\bf False} (the {\it anti-tautology}) corresponding to the {\it empty set}.

In addition to the above constant Boolean functions, there are three {\it atomic} functions.
The simplest is the {\it unary} function  {\bf NOT}, denoted by $\bar{x}$, that is defined by
$$
\bar{x} \, = \, \cases{ False \, , \quad if \quad x=True; \cr
                    True \, , \quad if \quad x=False \, .}
$$

The two other fundamental Boolean functions are the (inclusive) {\bf OR}, denoted by $\vee$ and 
{\bf AND}, denoted by $\wedge$. $x \vee y$ is True  unless both $x$ and $y$ are false, and
$x \wedge y$ is true only when both $x$ and $y$ are true.

By iterating these three operations on $n$ variables, one can get many {\it Boolean expressions}, and each
Boolean function has many possible expressions. 

From now on we will denote, as usual, {\it true} by $1$ and {\it false} by $0$.  Also let $x^1=x$ and $x^0=\bar{x}=1-x$.

One particularly simple type of expression is a {\bf conjunction} (also called {\it term}).
It is anything of the form, for some $t$, called its {\bf size},
$$
x_{i_1}^{j_1} \wedge \cdots  \wedge x_{i_t}^{j_t} \quad,
$$
where $1 \leq i_1<\dots<i_t \leq n$ and $j_i \in \{0,1\}$ for all $1 \leq i \leq t$.

Of interest to us in this article is the type of expression  called the {\it Disjunctive Normal Form} (DNF)
(featured prominently, along with its {\it dual}, {\it Conjunctive Normal Form}, (CNF),
in Norbert Blum's brave attempt ([Bl], see also [Zeil2])). It  simply  has the form
$$
\bigvee_{i=1}^{N} C_i \quad ,
$$
where each $C_i$ are pure conjunctions. 

Every Boolean expression corresponds to a unique function, but every function can be expressed in many ways,
and even in many ways that are DNF. One way that is the most straightforward way is the {\bf canonical DNF} form
$$
\bigvee_{\{ v \in f^{-1} (1) \}} \, \bigwedge_{i=1}^{n} \,  x_i^{v_i} \quad .
$$

Note that a pure conjunction of length $t$
$$
x_{i_1}^{j_1} \wedge \cdots  \wedge x_{i_t}^{j_t} \quad
$$
corresponds to a {\bf sub-cube} of dimension $n-t$, namely to
$$
\{\, (x_1, \dots, x_n) \, \vert \, x_{i_1}=j_1, \dots, x_{i_t}=j_t \} \quad .
$$

Hence, one can view a DNF as a (usually not exact) {\bf covering} of the set $f^{-1}(1)$ of truth-vectors 
by sub-cubes. In particular, a {\bf DNF tautology} is a covering of the whole $n$-dimensional unit
cube by lower-dimensional sub-cubes.

{\bf Digression: DNFs and the Million Dollar Problem}

The most fundamental problem in theoretical computer science, the question of whether {\bf P} is {\bf not} {\bf NP}
(of course it is not, but proving it rigorously is another matter), is equivalent to the question of whether there exists a polynomial time algorithm that decides if a given 
{\bf Disjunctive Normal Form} expression is the {\bf tautology} (i.e. the constant function $1$). 
Of course, there is an obvious {\it algorithm}:
For each term, find the truth-vectors covered by it, take the union, and see whether it contains all the
$2^n$ members of $\{0,1\}^n$. But this takes {\bf exponential time} and {\bf exponential memory}.

\vfill\eject

{\bf The Covering System Analog}

Input a system of congruences
$$
a_i \, (\, mod \, m_i \, ) \quad , \quad 1 \leq i \leq N \quad,
$$
and decide, in {\bf polynomial time}, whether it is a covering system. Initially it seems that we need
to check infinitely many cases, but of course (as already noted above), it suffices to check whether every integer between $1$ and
$lcm\,(\, m_1, \dots , m_N)$ belongs to at least one of the congruences. This seems fast enough!
Alas, the size of the input is the sum of the number of digits of the $a_i$'s and $m_i$'s and this
is less than a constant times the {\it logarithm} of $lcm\,(\, m_1, \dots , m_N)$, so just like for Boolean functions,
the naive algorithm is exponential time (and space) in the {\it input size}.

We next consider Boolean function analogs of covering systems.
The first one to consider such analogs was Melkamu Zeleke ([Zel]). Here we continue his pioneering work.

{\bf Boolean Function Analogs of Covering Systems}

We saw that a DNF tautology is nothing but a covering of the $n$-dimensional unit cube $\{0,1\}^n$ by
sub-cubes. So it is the analog of a covering system.

The analog of {\bf exact} covering systems is obvious: all the terms should cover disjoint sub-cubes. For example,
when $n=2$,  (from now on $xy$ means $x \wedge y$)
$$
x_1 x_2 \, \vee \, x_1 \bar{x_2} \, \vee \, \, \bar{x_1} {x_2} \, \vee \,  \bar{x_1} \bar{x_2} \quad,
$$
$$
x_1 \vee \, \, \bar{x_1} {x_2} \, \vee \,  \bar{x_1} \bar{x_2} \quad,
$$
are such.

In order  to define {\bf distinct} DNF, we define the {\bf support} of a conjunction as the set of the variables that
participate. For example , the support of the term $\bar{x_1} \bar{x_3} x_4 x_6$ is the set $\{x_1,x_3,x_4,x_6\}$.
In other words, we ignore the negations. For each $t$-subset of $\{x_1, \dots, x_n\}$ there are $2^t$ conjunctions with
that support. Geometrically speaking, two terms with the same support correspond to sub-cubes which are ``parallel'' to each other.

Note that the supports correspond to the modulo, $m$, and the assignments of negations (or no negation)
corresponds to a residue class modulo $m$.

A DNF tautology is {\bf distinct} if it has distinct supports.

An obvious example of a distinct DNF tautology in $n$ variables is
$$
 \bigvee_{i=1}^{n} \, x_i  \,\, \vee \, \, \wedge_{i=1}^{n} \bar{x_i}  \quad .
$$

More generally, for every $1 \leq t \leq n$, ($t \neq n/2$) the following is a distinct DNF tautology:

$$
\left ( \, \bigvee_{1 \leq i_1< i_2 < \dots <i_t \leq n} 
\, x_{i_1} \cdots x_{i_t} \, \right ) \, \vee \, 
\left ( \, \bigvee_{1 \leq j_1< j_2 < \dots <j_{n-t} \leq n} 
\, \bar{x}_{j_1} \cdots \bar{x}_{j_{n-t}} \, \right ) \quad .
$$
This follows from the fact that by the pigeon-hole principle, every $0-1$ vector of length $n$ either has
$\geq t$ $1$'s or $\geq n-t$ $0$'s.

The Boolean analog of the Mirsky-Newman-Davenport-Rado theorem is almost trivial. First, suppose we have an exact DNF tautology where the largest support has size $n$. That corresponds to a point (a $0$-dimensional subcube).
If it is the only one, then since a conjunction of length $t$ covers $2^{n-t}$ points, if all the other ones
are strictly smaller than $n$, and since they are all disjoint, they cover an even number of points, hence
there is no way that an exact DNF tautology would only have one term of size $n$.

If the largest size of a term is $<n$, then by projecting on appropriate sub-boxes one can reduce it
to the former case, and see that it must have a mate.

{\bf The Boolean Analog of  the Erd\H{o}s problem  is obviously TRUE}

Taking $n$ to be odd, the above DNF tautology with $t=(n-1)/2$  has ``minimal moduli'' (supports) of size $(n-1)/2$, and
that can  be made as large as one wishes.

{\bf First Challenge }

This leads to a more challenging problem: For each specific $n$, how large can the minimum clause size, let's call it $k$, in a distinct DNF tautology, be? 

An obvious {\bf necessary condition}, on {\bf density grounds}, is that
$$
\sum_{i=k}^{n} {{n} \choose {i}} \, \frac{1}{2^i} \, \geq 1  \quad .
$$

(Each subset of size $i$ of $\{1, \dots , n\}$ can only show up once and covers $2^{n-i}$ vertices of
the $n$-dimensional unit cube. Now use Boole's inequality that says that the number of elements of a union of
sets is $\leq$ than the sum of their cardinalities).

Let $A_n$ be the largest such $k$. The first $14$ values of $A_n$ are
$$
1, 1, 1, 2, 3, 4, 4, 5, 6, 7, 7, 8, 9, 10 \quad .
$$

We were able to find such optimal distinct DNF  tautologies for all $n\leq 14$ except for $n=10$, where the best
that we came up with was one that covers $1008$ out of the $1024$ vertices of the $10$-dimensional unit cube,
leaving $16$ points uncovered, and for $n=14$, where $276$ out of the $2^{14}=16384$ points were left uncovered.

See the output file 

{\tt http://sites.math.rutgers.edu/\~{}zeilberg/tokhniot/odt2.txt} \quad .

{\bf Second Challenge }

Another challenge is to come up with distinct DNF tautologies with all the terms of the {\bf same} size.
By density arguments a necessary condition for the existence of such a distinct DNF tautology
$$
 {{n} \choose {m}} \, \frac{1}{2^m} \geq 1 \quad .
$$
Let $B_m$ be the largest such $m$. The first $14$ values are
$$
 0, 0, 1, 2, 3, 3, 4, 5, 6, 6, 7, 8, 9, 9 \quad .
$$

Obviously for $n=3$, where $B_3=1$, it is not possible, since $x_1 \vee x_2 \vee x_3$ can't cover everything.
We were also unable to find such optimal DNF tautologies for $n=5$, where $B_5=3$ and we had to leave
one vertex uncovered, $n=9$, (with $B_9=6) $, where $13$ vertices were left uncovered, and $n=13$
(with $B_{13}=9$) where $2^{13}-8090=102$ vertices were left uncovered. For the other cases with $n\leq 14$, we met the
challenge. See the output file

{\tt http://sites.math.rutgers.edu/\~{}zeilberg/tokhniot/odt1.txt} \quad .

{\bf Supporting Maple Packages and Output}

Many more examples can be gotten from the Maple package

$\bullet$ {\tt http://www.math.rutgers.edu/\~{}zeilberg/tokhniot/dt.txt} \quad ,

whose output files are available from the {\it front} of this article

{\tt http://www.math.rutgers.edu/\~{}zeilberg/mamarim/mamarimhtml/dt.html} \quad .

{\bf The General Problem: Covering a Discrete Box by Non-Parallel Sub-boxes}

Let $\{a_i\}_{i=1}^{\infty}$ be a weakly increasing sequence of positive integers, with $a_1 \geq 2$.

Is it true that for every $m$ there exists an $n$ such that the box $[1,a_1] \times \dots [1,a_n]$ can be covered by
{\bf non-parallel} sub-boxes, each of dimension $\leq n-m$?.

We saw that for the Boolean case, with $a_i=2$ for each $i$ (and analogously, for each constant sequence), the answer
is obviously {\bf yes}.

On the other hand, if
$$
\sum_{i=1}^{\infty} \frac{1}{a_i} \, < \, \infty \quad ,
$$
the answer is obviously {\bf no}, since
$$
\prod_{i=1}^{\infty} (1+ \frac{1}{a_i}) \, < \, \infty \quad ,
$$
and by a trivial density argument, all tails of the product will eventually be less than $1$, so there is not enough room.

Regarding the original  Erd\H{o}s problem, Hough ([H]) proved the answer is {\bf no}  in the case with $a_i=p_i$, the
sequence of prime numbers.
(In fact, Hough proved the slightly harder result where the  moduli are not necessarily square-free.)
Here the sum of the reciprocals
{\it almost} converges. The very naive  Boole's inequality does not suffice to rule  out
a positive answer to the Erd\H{o}s problem,
but the Lov\'asz Local Lemma [that is also fairly weak; for example, it barely improves
the lower bounds for the Ramsey numbers] suffices to do the job.

So prime numbers were indeed {\bf red herrings}. All that was needed was their asymptotic behavior.
It would be interesting to see to what extent  Hough's  proof of impossibility extends to other
sequences $(a_i)$ for which the answer is neither an obvious Yes, nor an obvious No.
We also desperately need more powerful sieves than Lov\'asz, Brun, Selberg, and the other known
sieves, which, with lots of ingenuity, got close to the twin-prime conjecture (Yitang Zhang);
but even the still open twin-prime conjecture, and the Goldbach conjecture, are much weaker than the true state of affairs.

Even though the Bonferroni sieve is fairly weak, it can often decide {\it satisfiability}
(both positively and negatively). See the article [Za] by the first-named author.

{\bf References}
 
[BFF1]  Marc A. Berger, Alexander  Felzenbaum and Aviezri Fraenkel,
{\it 
   A nonanalytic proof of the Newman-Znam result for disjoint
   covering systems},   Combinatorica {\bf 6} (1986), 235-243.
 
[BFF2]  Marc A. Berger, Alexander  Felzenbaum and Aviezri Fraenkel,
{\it New results for covering systems of residue sets},
   Bull. Amer. Math. Soc. (N.S.) {\bf 14} (1986), 121--126.

[Bl] Norbert Blum, {\it A solution of the P versus NP problem}, \hfill\break
{\tt https://arxiv.org/abs/1708.03486} \quad .

[Bo] George Boole, L.L.D.,  ``{\it Investigations of 
THE LAWS OF THOUGHT, Of Which Are Founded
The Mathematical Theories of Logic and Probabilities}'', 
Macmillan, 1854. Reprinted by Dover, 1958.

[D] Jean-Paul Delahaye, {\it Cinq \'enigmes pour la rentr\'ee}, Logique et Calcul column, {\it Pour La Science},
No. {\bf 479} (Sept. 2017), 80-85.

[E1] Paul  Erd\H{o}s, {\it On integers of the form $2^k+p$ and some related problems}, Summa Brasil. Math. {\bf 2} (1950), 113-123.

[E2] Paul  Erd\H{o}s, {\it On a problem concerning covering systems}
(Hungarian, English summary), Mat. Lapok {\bf 3} (1952), 122-128.

[H] Bob Hough, {\it Solution of the minimum modulus problem for covering systems},  Ann. Math. {\bf 181} (2015), 361-382. \quad
{\tt https://arxiv.org/abs/1307.0874} \quad .

[S] Jamie Simpson, {\it Exact covering of the integers by arithmetic
progressions},
Discrete Math.  {\bf 59} (1986), 181-190.
 
[SZel] Jamie Simpson and Melkamu Zeleke, {\it On disjoint covering systems
with exactly one repeated modulus}, 
Adv. Appl. Math. {\bf 23} (1999), 322-332.

[W] Peter Winkler, ``{\it Mathematical Mind-Benders},'' A.K. Peters/CRC Press, 2007. 

[Za]  Anthony Zaleski, {\it Solving satisfiability using inclusion-exclusion}, \hfill\break
{\tt http://sites.math.rutgers.edu/\~{}az202/Z/sat/sat.pdf} . Maple package:  \hfill\break
{\tt  http://sites.math.rutgers.edu/\~{}az202/Z/sat/sat.txt } \quad .

[Zei1] Doron Zeilberger, {\it How Berger, Felzenbaum and Fraenkel revolutionized Covering Systems the same way
that George Boole revolutionized Logic}, Elect. J. Combinatorics {\bf8(2)} (2001) (special issue in honor of Aviezri Fraenkel), A1. \hfill\break
{\tt http://sites.math.rutgers.edu/\~{}zeilberg/mamarim/mamarimhtml/af.html} \quad .

[Zei2] Doron Zeilberger, {\it CNF-DNF and all that (Videotaped Lecture)}, available from \hfill\break
{\tt http://sites.math.rutgers.edu/\~{}zeilberg/mamarim/mamarimhtml/CNFDNFLecture.html} \quad .

[Zel] Melkamu Zeleke, {\it Ph.D. dissertation}, Temple University, 1998.

\bigskip
\hrule
\bigskip
Anthony Zaleski, Department of Mathematics, Rutgers University (New Brunswick), Hill Center-Busch Campus, 110 Frelinghuysen
Rd., Piscataway, NJ 08854-8019, USA. \hfill\break
Email: {\tt az202 at math dot rutgers dot edu}   \quad .
\bigskip
Doron Zeilberger, Department of Mathematics, Rutgers University (New Brunswick), Hill Center-Busch Campus, 110 Frelinghuysen
Rd., Piscataway, NJ 08854-8019, USA. \hfill\break
Email: {\tt DoronZeil at gmail dot com}   \quad .
\bigskip
{\bf Jan. 15, 2018}

\end